\magnification=1200


\hsize=125mm             
\vsize=195mm             
\parskip=0pt plus 1pt    
\clubpenalty=10000       
\widowpenalty=10000      
\frenchspacing           
\parindent=8mm           
			 
\let\txtf=\textfont
\let\scrf=\scriptfont
\let\sscf=\scriptscriptfont

\font\frtnrm =cmr12 at 14pt
\font\tenrm  =cmr10
\font\ninerm =cmr9
\font\sevenrm=cmr7
\font\fiverm =cmr5

\txtf0=\tenrm
\scrf0=\sevenrm
\sscf0=\fiverm

\def\rm{\fam0 \tenrm}

\font\frtnmi =cmmi12 at 14pt
\font\tenmi  =cmmi10
\font\ninemi =cmmi9
\font\sevenmi=cmmi7
\font\fivemi =cmmi5

\txtf1=\tenmi
\scrf1=\sevenmi
\sscf1=\fivemi

 \def\oldstyle{\fam1 \tenmi}

\font\tensy  =cmsy10
\font\ninesy =cmsy9
\font\sevensy=cmsy7
\font\fivesy =cmsy5

\txtf2=\tensy
\scrf2=\sevensy
\sscf2=\fivesy

\font\tenit  =cmti10
\font\nineit =cmti9
\font\sevenit=cmti7
\font\fiveit =cmti7 at 5pt

\txtf\itfam=\tenit
\scrf\itfam=\sevenit
\sscf\itfam=\fiveit

\def\it{\fam\itfam\tenit}

\font\tenbf  =cmb10
\font\ninebf =cmb10 at  9pt
\font\sevenbf=cmb10 at  7pt
\font\fivebf =cmb10 at  5pt

\txtf\bffam=\tenbf
\scrf\bffam=\sevenbf
\sscf\bffam=\fivebf

\def\bf{\fam\bffam\tenbf}

\newfam\msbfam       
\font\tenmsb  =msbm10
\font\sevenmsb=msbm7
\font\fivemsb =msbm5

\txtf\msbfam=\tenmsb
\scrf\msbfam=\sevenmsb
\sscf\msbfam=\fivemsb

\def\msb{\fam\msbfam\tenmsb}
\def\Bbb#1{{\msb #1}}

\newfam\scfam
\font\tensc  =cmcsc10

\txtf\scfam=\tensc

\def\sc{\fam\scfam\tensc}


\def\frtnmath{%
\txtf0=\frtnrm        
\txtf1=\frtnmi         
}

\def\frtnpoint{%
\baselineskip=16.8pt plus.5pt minus.5pt%
\def\rm{\fam0 \frtnrm}%
\def\oldstyle{\fam1 \frtnmi}%
\everymath{\frtnmath}%
\everyhbox{\frtnrm}%
\frtnrm }


\def\ninemath{%
\txtf0=\ninerm        
\txtf1=\ninemi        
\txtf2=\ninesy        
\txtf\itfam=\nineit      
\txtf\bffam=\ninebf      
}

\def\ninepoint{%
\baselineskip=10.8pt plus.1pt minus.1pt%
\def\rm{\fam0 \ninerm}%
\def\oldstyle{\fam1 \ninemi}%
\def\it{\fam\itfam\nineit}%
\def\bf{\fam\bffam\ninebf}%
\everymath{\ninemath}%
\everyhbox{\ninerm}%
\ninerm }



\def\text#1{\hbox{\rm #1}}

\def\cite#1{{\uppercase{#1}}}
\def\ref#1{{\uppercase{#1}}}
\def\label#1{{\uppercase{#1}}}

\def\br{\hfill\break} 



\def\topmatter{\null\firstpagetrue\vskip\bigskipamount}
\def\endtopmatter{\vskip2\bigskipamount}

\def\title#1{%
\vbox{\raggedright\frtnpoint
\noindent #1\par}
\vskip 2\bigskipamount}        


\def\shorttitle#1{\rightheadtext={#1}}              

\newif\ifThanks
\global\Thanksfalse

\def\author#1{\begingroup\raggedright
\noindent{\sc #1\ifThanks$^*$\else\fi}\endgroup
\leftheadtext={#1}\vskip \bigskipamount}

\def\address#1{\begingroup\raggedright
\noindent{#1}\endgroup\vskip\bigskipamount}

\def\endabstract{\endgroup}

\long\def\abstract#1\endabstract{\par
\begingroup\ninepoint\narrower
\noindent{\sc Abstract.\enspace}#1%
\vskip\bigskipamount\endabstract}

\def\section#1#2{\bigbreak\bigskip\begingroup\raggedright
\noindent{\bf #1.\quad #2}\nobreak
\medskip\endgroup\noindent\ignorespaces}

\def\proclaim#1{\medbreak\noindent{\sc #1.\enspace}\begingroup
\it\ignorespaces}
\def\endproclaim{\endgroup\bigbreak}

\def\remark#1{\medbreak\noindent{\sc Remark \enspace}
\begingroup\ignorespaces}
\def\endremark{\endgroup\bigbreak}


\def\qed{$\mathord{\vbox{\hrule\hbox{\vrule
		               \hskip5pt\vrule height5pt\vrule}\hrule}}$}

\def\demo#1{\medbreak\noindent{{#1}.\enspace}\ignorespaces}
\def\enddemo{\penalty-100\null\hfill\qed\bigbreak}

\newdimen\EZ

\EZ=.5\parindent

\newbox\itembox

\newdimen\ITEM
\newdimen\ITEMORG
\newdimen\ITEMX
\newdimen\BUEXE

\def\iteml#1#2#3{\par\ITEM=#2\EZ\ITEMX=#1\EZ\BUEXE=\ITEM
\advance\BUEXE by-\ITEMX\hangindent\ITEM
\noindent\leavevmode\hskip\ITEM\llap{\hbox to\BUEXE{#3\hfil$\,$}}%
\ignorespaces}


\newif\iffirstpage\newtoks\righthead
\newtoks\lefthead
\newtoks\rightheadtext
\newtoks\leftheadtext
\righthead={\ninepoint\rm\hfill{\the\rightheadtext}\hfill\llap{\folio}}
\lefthead={\ninepoint\rm\rlap{\folio}\hfill{\the\leftheadtext}\hfill}
\headline={\iffirstpage\hfill\else
\ifodd\pageno\the\righthead\else\the\lefthead\fi\fi}
\footline={\iffirstpage\hfill\global\firstpagefalse\else\hfill\fi}

\leftheadtext={}
\rightheadtext={}


\def\Refs{\bigbreak\bigskip\noindent{\bf References}\medskip
\begingroup\ninepoint\parindent=40pt}
\def\endRefs{\par\endgroup}
\def\endref{}

\def\ref{\par}
\def\key#1{\item{\hbox to 30pt{[#1]\hfill}}}
\def\by{}
\def\paper{\it}
\def\jour{\rm}
\def\vol{\bf}
\def\pages{\rm}
\def\yr{}
\def\publaddr{\rm}
\def\publ{}


\def\cline#1{\leftline{\hfill#1\hfill}}

\def\bR{\Bbb R}
\def\bC{\Bbb C}

\def\bZ{\Bbb Z}

\input epsf

\topmatter
\title{Real deformations and complex topology \br
of plane curve singularities}
\author{Norbert A'Campo}
\shorttitle{Real deformations and complex topology.}
\endtopmatter

\par
\noindent
{\bf Table of contents}\par\noindent
1. Introduction \br
2. Real deformations of plane curve singularities \br
3. Complex topology of plane curve singularities \br
4. The singularity $D_5$ and a graphical algorithm in general
\br
5. An example of global geometric monodromy \br
6. Connected divides and fibered knots. Proof of Theorem 2 \br

\section{1}{Introduction}
The geometric monodromy $T$ of a  curve singularity in the complex
plane is a
diffeomorphism of a compact surface with boundary $(F,\partial F)$
inducing
the identity on the boundary, which is well defined up to isotopy
relative to
the boundary. The geometric monodromy of a  curve singularity in the
complex
plane determines the local topology of the singularity.  As element of
the
mapping class group of the surface $(F,\partial F)$, the diffeomorphism
$T$
can be written as a composition of Dehn twists. In section 3 of this
paper the
geometric monodromy of an isolated plane curve singularity is written
explicitly as a
composition  of right Dehn twists. In fact, a global graphical
algorithm
for the construction of the surface $(F,\partial F)$ with a system of
simply closed curves on it is given in section 4, such that the curves
of this system are the  vanishing cycles of a real morsification of the
singularity. In section 5, as an illustration, the global geometric
monodromy
of the polynomial $y^4-2y^2x^3+x^6-x^7-4yx^5,$ which has two critical
fibers, is computed. \br
The germ of a curve singularity in $\bC^2$ is a finite union of
parametrized
local branches $b_i:\bC \to \bC^2, 1 \leq i \leq r.$ First observe,
that
without loss of generality for the local topology, we can assume that
the
branches have a  real polynomial parameterization. The combinatorial
data used
to describe the geometric monodromy of a curve singularity come from
generic
real polynomial deformations of the parameterizations of the local
branches
$b_{i,t}:\bC \to \bC^2, 1 \leq i \leq r, t \in [0,1]$, such
that:\br
(i) $b_{i,0}=b_i, 1 \leq i \leq r,$\br
(ii) for some $\rho > 0$ the intersection of the union of the branches
with the $\rho-$ball $B$ at the singular point of curve in $\bC^2$ is
a representative of the germ of the curve and $B$ is a Milnor ball for
the
germ,\br
(iii) the images of $b_{i,t},1 \leq i \leq r, t \in [0,1],$ intersect
the
boundary of the ball $B$ transversally,\br
(iv) the union of the images $b_{i,t}(\bR), 1 \leq i \leq r,$ has for
every $t \in (0,1]$ the maximal possible number of double points in the
interior of $B$.\br

Such deformations correspond to real morsifications of the defining
equation
of the singularity and were used to study the local monodromy
in  [AC2],[AC3],[G-Z]. Real deformations of singu\-larities of plane
algebraic curves with the maximal possible number of double points in
the real plane were
discovered by Charlotte Angas Scott [S1,S2]. I thank Egbert Brieskorn
for
having drawn my attention on the references [S1],[S2].

In section 6 we start with  a connected divide, which  defines as
explained in section 3 a classical link. We will construct a map from
the complement of the link of a connected divide to the circle and
prove that this map is a fibration.
This fibration is for a divide  of a plane curve singularity a model
for the Milnor fibration of the singularity. The link of most connected
divides  are  hyperbolic. In a forth coming paper we will study the
geometry of a link of a divide.

We used MAPLE for the drawings of parametrized curves and for the
computation of suitable deformations of the polynomial equations.
Of great help for the investigation of topological changes in
families of polynomial equations is the mathematical software SURF
which
has been developed by Stefan Endrass. I thank  Stefan Endrass warmly
for
permitting me to use SURF. Part of this work was done in Toulouse and
I thank the members of the Laboratory \'Emile Picard for their
hospitality.

\par\noindent

\section{2}{\bf Real deformations of plane curve singularities}

Let $f:\bC^2 \to \bC$ be the germ at $0 \in \bC^2$ of an holomorphic
map
with $f(0)=0$ and having an isolated singularity $S$ at $0$. We are
mainly
interested in the study of topological properties of  singularities,
therefore we can assume without  loss of generality that the germ $f$
is a
product of locally irreducible real polynomials. Having chosen a Milnor
ball
$B(0,\rho)$ for $f$, there exists a real polynomial deformation family
$f_t,t \in [0,1],$ of $f$ such that for all $t$ the $0$-level of $f_t$
is transversal to the boundary of the ball $B(0,\rho)$ and such that
for
all $t \in ]0,1]$ the $0$-level of $f_t$ has  $\delta$ transversal
double
points in the interior of the disk  $D(0,\rho):=B(0,\rho) \cap \bR^2$,
where
the Milnor number $\mu$ and the number $r$ of local branches of $f$
satisfy
$\mu=2\delta-r+1.$ In particular, the $0$-level of $f_t, t \in ]0,1],$
has
in  $D(0,\rho)$ no self tangencies or triple intersections. It is
possible
to choose for $f_t,t \in [0,1],$ a family
of defining equations for the union of the images of
$b_{i,t}, 1 \leq i \leq r.$ The deformation $f_t,t \in [0,1]$ is called
a
real morsification with respect to the Milnor ball $B(0,\rho)$ of $f.$
So, the $0$-level of the restriction of
$f_t,t \in ]0,1],$ to $D(0,\rho)$ is an immersion without
self-tangencies and
having only transversal self-intersections of $r$ copies of an interval
(see [AC2],[AC3],[G-Z]). The $0$-level of the restriction of
$f_t,t \in ]0,1],$ to $D(0,\rho)$ is up to a diffeomorphism independent
of
$t,$ it is called a divide ("partage" in [AC2]) and it is shown that
for
instance the divide determines the homological monodromy group of the
versal deformation of the singularity. Figure 1 represents a divide for
the singularity at $0 \in \bC^2$ of the curve $(y^5-x^3)(x^5-y^3)=0$.
\midinsert
\cline{\epsffile{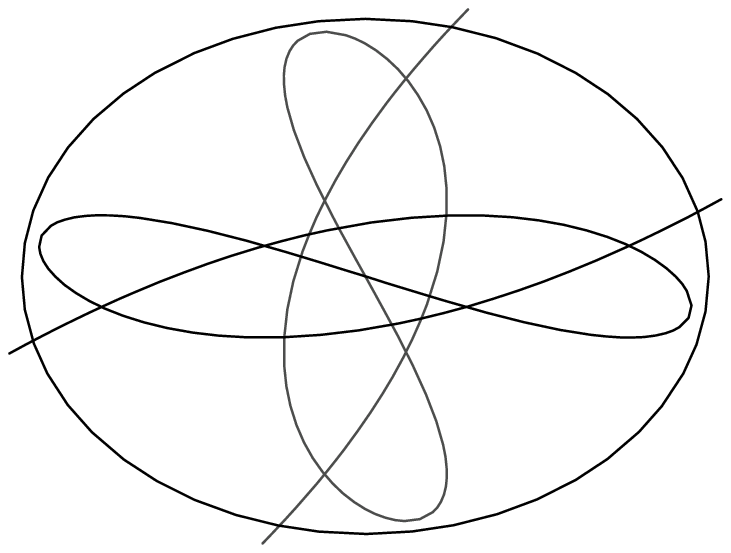}}
\medskip
\centerline{Figure 1: A divide for $(y^5-x^3)(x^5-y^3)=0.$}
\endinsert

\remark{\bf Remark}
The transversal isotopy class of the divide of a singularity with real
branches is not a topological invariant of the singularity. The
singularities of $y^4-2y^2x^3+x^6+x^7$ and $y^4-2y^2x^3+x^6-x^7$ have
congruent but not transversal isotopic divides. The singularities
$(x^2-y^2)(x^2-y^3)(y^2-x^3)$ and $(x^2-3xy+2y^2)(x^2-y^3)(y^2-x^3)$
are topologically equivalent but can not have congruent divides.
The singularity $y^3-x^5$ admits two divides, which give a model for
the smallest possible transition, according to the mod 4 congruence of
V. Arnold [A] and its 
celebrated strengthening to a mod 8 congruence  of 
V. A. Rohlin [R1,R2], of odd ovals to
even ovals for projective real M-curves of even degree. I owe this remark to
Oleg Viro [V].
More precisely, there exist  a polynomial family $f_s(x,y), s \in \bR,$ of
polynomials of degree $6,$ having the central symmetry 
$f_s(x,y)=-f_{-s}(-x,-y)$ such that the levels
$f_s(x,y)=0,\ s \not=0,$
are divides for the singularity $f_0(x,y)=y^3-x^5.$ Moreover, for $s \in \bR,\
s\not=0,$
the divide $f_s(x,y)=0$ has for regions on which the function has the sign of
the parameter $s.$ Therefor,
at
$s=0$ four regions of $f_s(x,y)=0$ collapse and hence, four ovals of
$f_s(x,y)=(s/2)^E$ collapse and change parity 
at $s=0$ if the exponent $E$ is big and odd. 

\midinsert
\cline{\epsffile{divide_links.eps} \ \epsffile{divide_rechts.eps}}
\medskip
\centerline{Figure 2a: The divides $f_{\pm 1}(x,y)=0.$}
\endinsert

In Figure $2b$ are drawn the
smoothings of the
divides $f_{-1}(x,y)=0$ and $f_{+1}(x,y)=0$ of 
the singularity of $f_0(x,y)=y^3-x^5$ at $0$. Each of the smoothings 
$\{f_{\pm1}(x,y)=\pm \epsilon \}\cap D$  consists of four ovals and a 
chord, such that the ovals lie
on the positive, respectively on the negative side, of the chord. 
Such a family $f_s(x,y)$
is for instance given by:

\midinsert
\cline{\epsffile{links.eps} \ \epsffile{rechts.eps}}
\medskip
\centerline{Figure 2b: Four ovals change parity.}
\endinsert

$$
y^3-x^5
-{{125}\over{8}}\,{s}^{3}{x}^{6}
+({{375}\over{64}}\,{s}^{6}
+{{245}\over{16}}\,{s}^{4}
-{{25}\over{4}}\,{s}^{2}){x}^{5}
+{{75}\over{4}}\,{s}^{2}{x}^{4}y
$$
$$
+({{2695}\over{128}}\,{s}^{7}
+{{21625}\over{256}}\,{s}^{9}
-{{35}\over{8}}\,{s}^{3}
-{{4847}\over{160}}\,{s}^{5}){x}^{4}
-(5\,s
+{{159}\over{4}}\,{s}^{3}
-{{75}\over{4}}\,{s}^{5}){x}^{3}y
$$
$$
-{{15}\over{2}}\,s{x}^{2}{y}^{2}
-({{2703}\over{500}}\,{s}^{6}
+{{1281}\over{32}}\,{s}^{8}
-{{29625}\over{512}}\,{s}^{12}
+{{2345}\over{128}}\,{s}^{10})x^3
$$
$$
-({{17625}\over{256}}\,{s}^{8}
+{{3583}\over{32}}\,{s}^{6}
+{{5793}\over{400}}\,{s}^{4}){x}^{2}y
-({{95}\over{2}}\,{s}^{4}
+{{53}\over{10}}\,{s}^{2})x{y}^{2}
$$
$$
-({{997}\over{4000}}\,{s}^{9}
+{{42875}\over{2048}}\,{s}^{15}
+{{4575}\over{256}}\,{s}^{13}
+{{857}\over{320}}\,{s}^{11}){x}^{2}
$$
$$
-({{177325}\over{2048}}\,{s}^{11}
+{{1803}\over{200}}\,{s}^{7}
+{{35441}\over{512}}\,{s}^{9})xy
-({{6395}\over{128}}\,{s}^{7}
+{{317}\over{40}}\,{s}^{5}){y}^{2}
$$
$$
+({{19871}\over{1280}}\,{s}^{14}
+{{10165}\over{1024}}\,{s}^{16}
-{{59125}\over{4096}}\,{s}^{18}
+{{4171}\over{2000}}\,{s}^{12})x
$$
$$
+({{51025}\over{4096}}\,{s}^{12}
+{{54223}\over{25600}}\,{s}^{10}
-{{153725}\over{16384}}\,{s}^{14})y
$$

\endremark

\remark{\bf Problem}
Classify up to transversal isotopy, i.e.  isotopy through  immersions
with
only transversal double and triple point crossings, the   divides for
an
isolated real plane curve singularity.
\endremark
\goodbreak

\section{3}{Complex topology of plane curve singularities}
In this section we wish to explain how one can read off from the divide
of
a plane curve singularity $S$ the local link $L$, the Milnor fiber and
the
geometric monodromy group of the singularity. In particular, we will
give
the geometric monodromy of the singularity explicitly as a product of
Dehn twists.

Let $P \subset D(0,\rho)$ be the divide of the singularity $f$. For a
tangent
vector $v \in TD(0,\rho)=D(0,\rho) \times\bR^2$ of $D$ at the point
$p\in D(0,\rho)$ let $J(v) \in \bC^2$ be the point $p+iv.$ The Milnor
ball
$B$ can be viewed as
$$
B(0,\rho)=\{J(v) \mid v \in T(D(0,\rho)) \ \text{ and }\ \|J(v)\| \leq
\rho\}.
$$
Observe that
$$L(P):=\{J(v) \mid v \in T(P) \ \text{ and }\  \|J(v)\|=\rho \}$$
is a closed  submanifold of dimension one in the boundary of the Milnor
ball $B(0,\rho)$. We call $L(P)$ the link of the divide $P$.
Note further  that
$$R(P):=\{J(v) \mid v \in T(P) \ \text{ and }\  \|J(v)\| \leq \rho \}$$
is an immersed surface  in $B(0,\rho)$ with boundary $L(P)$  having
only
transversal double point singularities.  Let $F(P)$ be the surface
obtained from $R(P)$ by replacing the local links of its  singularities
by cylinders. The differential model of those replacements is as
follows:
let $\chi: \bC^2 \to \bR$ be a smooth bump function at $0 \in \bC^2$;
replace the immersed surface
$\{(x,y) \in \bC^2 \mid xy=0\}$ by the smooth surface
$\{(x,y) \in \bC^2 \mid xy=\tau^2\chi(x/\tau,y/\tau)\}$, where $\tau$
is a
sufficiently small positive real number. We call $R(P)$ the singular
and
$F(P)$ the regular ribbon surface of the divide $P.$ The connected,
compact
surface $F(P)$ has genus $g:=\delta-r+1$ and $r$ boundary components.
Note, that $g$ is the number of regions of the divide $P$. A region of
$P$
is a  connected component of the complement of $P$ in $D(0,\rho)$,
which
lies in the interior of $D(0,\rho)$. For the example drawn in Figure
$1,$ we
have $r=2,\delta=17,g=16.$

The ribbon surface $R(P)$ carries a natural orientation, since
parametrized by an open subset of the tangent space $T(\bR).$ Hence the
surface $F(P)$ and the link $L(P)$ are also naturally oriented. We
orient $B$ as a submanifold of $-T\bR^2,$ which is the orientation of
$B$ as a submanifold in $\bC^2.$

\proclaim{Theorem 1}
Let $P$ be the divide for an isolated plane curve singularity $S.$  The
submanifold $(F(P),L(P))$ is up to isotopy a model for the Milnor
fiber of the singularity $S$.
\endproclaim

\demo{ \bf Proof} Choose $0 < \rho_{-} < \rho$ such that $P \cap
D(0,\rho_{-})$ is
still a divide for the singularity $S.$ Along the divide the singular
level
$F_{t,0}:=\{(x,y) \in B \mid f_t(x,y)=0 \}$ is up to order $1$ tangent
to the
immersed surface $R(P)$. Hence, for
$B_{-}':=\{u+iv \in B(0,\rho_{-}) \mid u,v \in
\bR^2, ||v|| \leq \rho'\}$ with $0 < \rho' << \rho,$ the intersections
$R'(P):=\partial B_{-}' \cap
R(P)$  and $F'_{t,0}:=\partial B_{-}' \cap F_{t,0}$ are
transversal and are regular
collar neighbourhoods of
the divide in $R(P)$ and in $F_{t,0}.$  Therefore the
nonsingular level $F_{t,\eta}:=\{(x,y) \in B(0,\rho) \mid f_t(x,y)=\eta
\}$,
where $\eta \in \bR$ is
sufficiently small, contains in its interior $F'_{t,\eta}:=B_{-}' \cap
F_{t,\eta},$ which is a diffeomorphic copy of the surface with
boundary $F(P)$. Since  $F'_{t,\eta}$ and $F_{t,\eta}$ are connected
surfaces both with $r$ boundary components and the
intersection forms on the first homology
are isomorphic, the difference $F_{t,\eta}\setminus F'_{t,\eta}$  is a
union of open collar
tubular neighbourhoods of the boundary components of the surface
$F_{t,\eta}$. So, the surfaces
$F_{t,\eta},F'_{t,\eta}$ and $F(P)$ are diffeomorphic.
We conclude by observing that the nonsingular levels $F_{t,\eta}$
and the Milnor fiber are connected in
the local unfolding  through nonsingular levels.
\enddemo

From this proof it follows also that the local link $L(S)$ of the
singularity
$S$ in
$\partial B$ is
cobordant  to the sub\-manifold $\partial F'_{t,0}$ in $\partial
B_{-}'.$ The
cobordism is given by the pair $(B \setminus
\text{int}(B_{-}'),F_{t,\eta}
\setminus \text{int}(F'_{t,\eta})).$ It
is clear, that the pairs $(\partial B,L(P)),$
$(\partial B_{-}',\partial F'_{t,0})$ and
$(\partial B_{-}',\partial F'_{t,\eta})$ are
diffeomorphic.
One can prove even more:

\proclaim{Theorem 2}
Let $P$ be the divide for an isolated plane curve singularity $S.$
The pairs $(\partial B,L(S))$, where $L(S)$ is the local link of the
singularity
$S,$ and $(\partial B,L(P))$ are diffeomorphic.
\endproclaim

The proof is given in section 6.

\remark{\bf Remark}
The signed planar Dynkin diagram of the divide determines up to isotopy
the
divide of  the singularity. It follows from Theorem 2, that the signed
planar
Dynkin diagram
determines geometrically the topology of the singularity. Using the
theorem
of Burau and Zariski stating that the topological type of a plane curve
singularity is determined by the mutual intersection numbers of the
branches
and the Alexander polynomial of each branch, the authors L. Balke and
R. Kaenders [B-K] have proved that the  signed   Dynkin diagram,
without its
planar embedding, determines the topology of the singularity.
\endremark

\midinsert
\cline{\epsffile{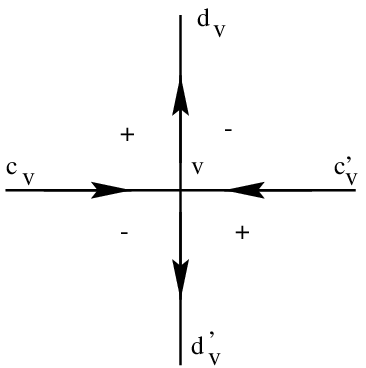}}
\medskip
\centerline{Figure $3:$ The link $P_v.$}
\endinsert

We need a combinatorial description of the surface $F(P).$ For a divide
$P$
we define: A vertex of $P$ is double point of $P$, and an edge of $P$
is the
closure of a connected component of the complement of the vertices in
$P$.
Now we choose an orientation of $\bR^2$, and a small deformation $\bar
f$
of the polynomial $f$ such that the $0$-level of $\bar f$ is the divide
$P.$
We call a region of the divide positive or negative according to the
sign
of $\bar f$. We orient the boundaries of the positive regions such that
the
outer normal and the oriented tangents of the boundary agree in this
order
with the chosen orientation of $\bR^2$. We choose a midpoint on each
edge,
which connects two vertices. The link $P_v$ of a vertex $v$ is the
closure
of the connected component of the complement of the midpoints in $P$
containing the given vertex $v$.

\midinsert
\cline{\epsffile{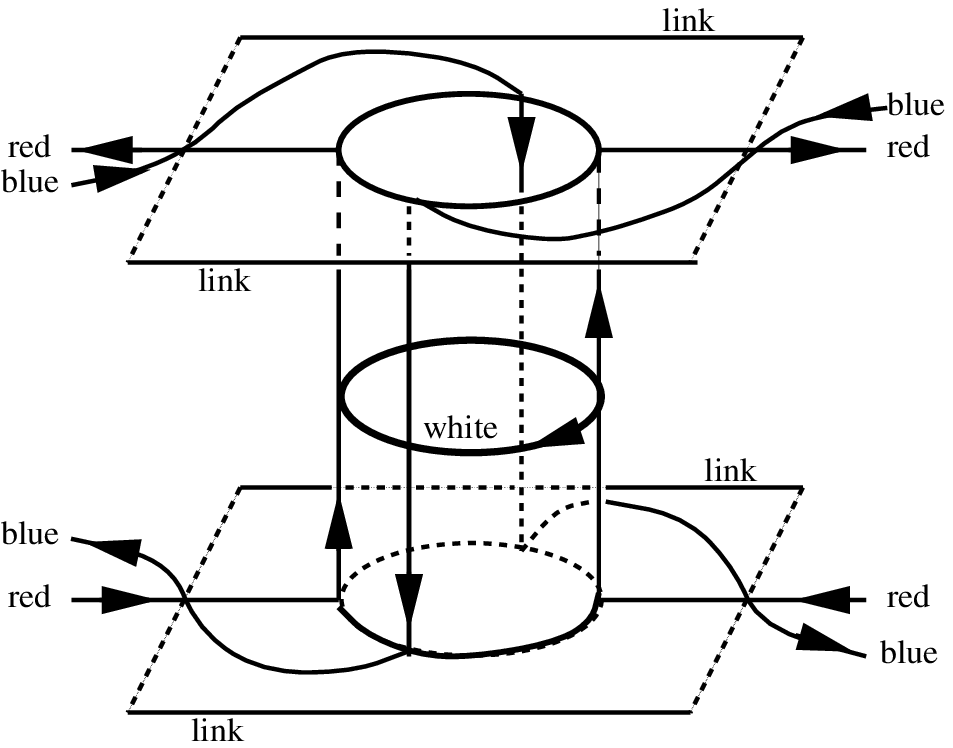}}
\medskip
\centerline{Figure $4:$ A piece of surface $F_v.$}
\endinsert

For each vertex $v$ of $P$ we will construct
a piece of  surface $F_v$, such that those pieces glue together and
build
$F(P)$. Let $P_v$ be the link of the vertex $v$. Call $c_v,c_v'$ the
endpoints of the branches of $P_v,$ which are oriented towards $v,$ and
$d_v,d_v'$ the endpoints of the branches of $P_v,$ which are oriented
away
from $v$. Thus, $c_v,c_v',d_v,d_v'$ are midpoints or endpoints of the
divide $P$ (see Figure $3.$) Using an orientation of the divide $P$ we
label
$c_v,c_v'$ such that $c_v'$
comes after $c_v$, and we label $d_v,d_v'$ such that the sector
$c_v,d_v$ is
in a positive region. Then $F_v$  is the surface with boundary and
corners
drawn in Figure $4.$ There are 8 corners and there are 8 boundary
components
in between the corners,
4 of them will get a marking by $c_v,c_v',d_v,d_v'$, which will
determine the gluing with the piece of the next vertex and 4 do not
have a marking. The
gluing of the pieces $F_v$ along the marked boundary components
according to
the gluing scheme given by the divide $P$ yields the surface $F(P)$. On
the
pieces $F_v$ we have drawn oriented curves colored red, white, and
blue.
The white  curves are simple closed pairwise disjoint curves. The
surface
$F(P)$ will be oriented such that the curves taken in the order
red-white-blue have nonnegative intersections.
The remaining red curves
glue together and build a red graph on  $F(P)$. The remaining blue
curves
build a blue graph. After deleting each contractible component of the
red or
blue graph, each of the remaining components contains a simple closed
red or
blue curve.   All together, we have constructed on $F(P)$ a system of
$\mu$
simple closed curves $\delta_1,\delta_2, ...
,\delta_{\mu-1},\delta_\mu$,
which we list by first taking red, then white and finally blue. We
denote
by $n_{+}$ the number of red curves which is also the number of
positive
regions, by $n_{\cdot}$ the number of crossing points and by $n_{-}$
the
number of blue which equals the number of negative regions of the
divide $P$.

Let $D_i$ be the right Dehn twist along the curve $\delta_i.$ A model
for
the right Dehn twist is  the linear action $(x,y) \mapsto (x+y,y)$ on
the
cylinder $\{(x,y) \in \bR/\bZ \times \bR \mid 0 \leq y \leq
1\}$ with as orientation the product of the natural orientations of the
factors. A right Dehn twist
around a simply
closed curve $\delta$ on an oriented surface is obtained by embedding
the
model as an oriented bicollar neighbourhood of $\delta$ such that
$\delta$
and $\bR/\bZ \times \{1/2\}$ of the model match.
The local geometric monodromy
of the singularity of $xy=0$ is as diffeomorphism a right Dehn
twist (voir le Th\'eor\`eme Fondamental, page 23 de [L], et page 95 de
[P-S]). Using as in [AC2] a local version of a  Theorem of Lefschetz,
one obtains:

\proclaim{Theorem 3}
Let $P$ be the divide for an isolated plane curve singularity $S.$ The
Dehn
twists $D_i$ are generators for the geometric monodromy group of the
unfolding
of the singularity $S$. The product $T:=D_{\mu}D_{\mu-1} ... D_2D_1$ is
the
local geometric monodromy of the singularity $S$.
\endproclaim

\section{4}{The singularity $D_5$ and a graphical algorithm in
 general}

We will work out the picture for the singularity $D_5$ with the
equation
$x(x^3-y^2)$ and the divide given by the deformation
$(x-s)(x^3+5sx^2-y^2), s \in [0,1],$ which is shown for $s=1$
in Figure $5.$ There are  one positive triangular region, one negative
region and three
crossings. By gluing three pieces together, one gets the Milnor fiber
with a system of
vanishing cycles as depicted in Figure $6.$

\midinsert
\cline{\epsffile{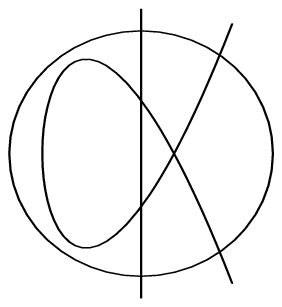}}
\medskip
\centerline{Figure $5:$  A divide for the singularity $D_5.$}
\endinsert

An easy and fast graphical algorithm of visualizing the Milnor fiber
with a
system of vanishing cycles directly from the divide is as follows:
think the divide
as  a road network which has $\delta$ junctions, and replace
every junction by a
roundabout,
which leads you to a new road network with $4\delta$ T-junctions.
Realize now
every road section in between two T-junctions by a strip with a half
twist. Do the same for every road section in between a
T-junction and the boundary of the divide. Altogether you will need
$6\delta+r$ strips. The core line of the four
strips of a roundabout is a white vanishing cycle, the strips
corresponding
to boundary edges and
corners of a positive or negative region have as core line a red or blue
vanishing cycle.

\midinsert
\cline{\epsffile{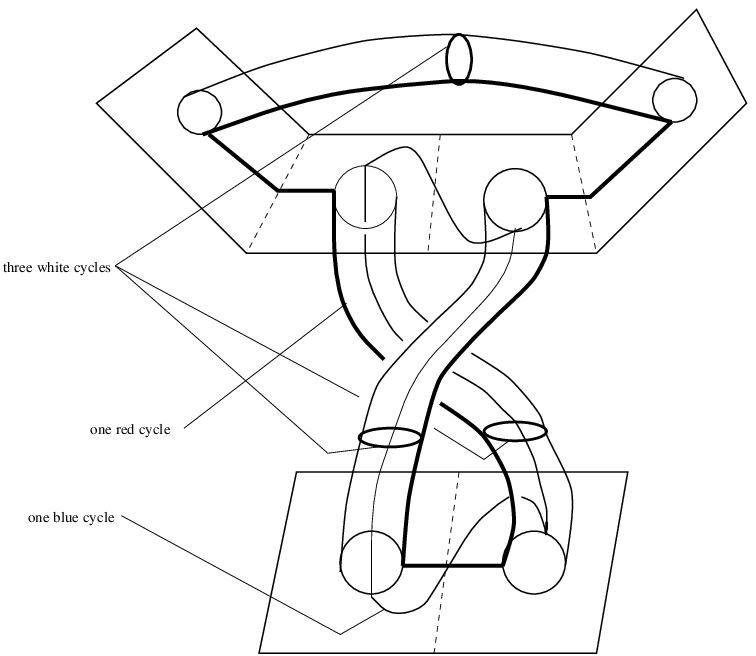}}
\medskip
\centerline{Figure $6:$ Milnor fiber with vanishing cycles for $D_5.$}
\endinsert

In Figure $7$ is worked out the singularity with two
Puiseux pairs and $\mu=16$, where we used the divide from Figure $9.$

\midinsert
\cline{\epsffile{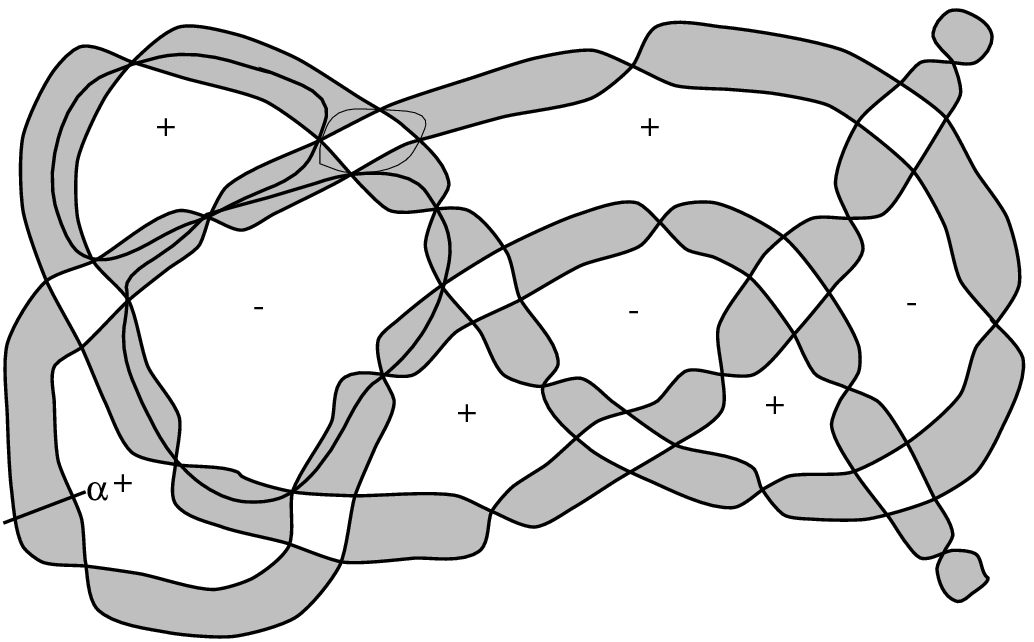}}
\medskip
\centerline{Figure $7:$ Milnor fiber with vanishing cycles for
$y^4-2y^2x^3+x^6-x^7-4yx^5.$}
\endinsert

We have drawn
for convenience in
Figure $7$ only one red, white, or blue cycle.
We have also indicated the position of the arc $\alpha,$ which will
play a role in the next section.
\goodbreak

\section{5}{An example of global geometric monodromy}
Let $b:\bC \to \bC^2, b(t):=(t^6+t^7,t^4)$ be the parametrized curve
$C$
having at $b(0)=(0,0)$ the singularity with two essential Puiseux pairs
and
with local link  the compound cable knot $(2,3)(2,3).$ The polynomial
$f(x,y):=y^4-2y^2x^3+x^6-x^7-4yx^5$ is the equation
of $C.$ The function $f:\bC^2 \to \bC$ has besides $0$ the only other
critical
value $c:=14130940973168155968/558545864083284007.$
The fiber of $0$ has besides its singularity at $(0,0)$ a nodal
singularity
at $(-8,-4)$, which corresponds to the node $b(-1+i)=b(-1-i)=(-8,-4).$
The
geometric monodromy of the singularity at $(0,0),$ which is up to
isotopy
piecewise of finite order, is described in [AC1]. The fiber of $c$ has
a
nodal singularity at $(1014/343,16807/79092).$ The singularity at
infinity
of the curve $C$ is at the point $(0:1:0)$ and its local equation is
$z^3-2z^2x^4+zx^6-x^7-4zx^5$, whose singularity is topologically
equivalent
to the singularity $u^3-v^7$ with Milnor number $12.$
The function $f$ has no critical values coming from infinity. We aim at
a
description of the global geometric monodromy of the function $f.$
Working with the distance on $\bC^2$ given by $||(x,y)||^2:= |x|^2+
4|y|^2,$
we have that the parametrized curve $b$ is transversal to the spheres
$S_r:=\{(x,y) \in \bC^2 \mid |x|^2+ 4|y|^2=r^2\}$ with center $0 \in
\bC^2$
and radius $r > 0.$ So for $0<r< 8\sqrt{2},$ the intersection
$K_r:=C \cap S_r$ is  the local knot in $S_r$ of the singularity at
$0 \in \bC^2$ (see Figure $8$), at $r=8\sqrt{2}$ the knot $K_r$ is
singular with one
transversal crossing, and for $8\sqrt{2} < r$ the knot $K_r$ is the so
called knot at infinity of the curve $C.$ The crossing at the bottom of Fig.
$8$ flips for $r=8\sqrt{2}$ and the knot $K_r, 8\sqrt{2} < r,$ 
becomes the $(4,7)$ torus knot. 
By making one extra total
twist
in a braid presentation of the knot $K_r$ one gets the local knot of
the
singularity at infinity of the projective completion of the curve $C.$

\midinsert
\cline{\epsffile{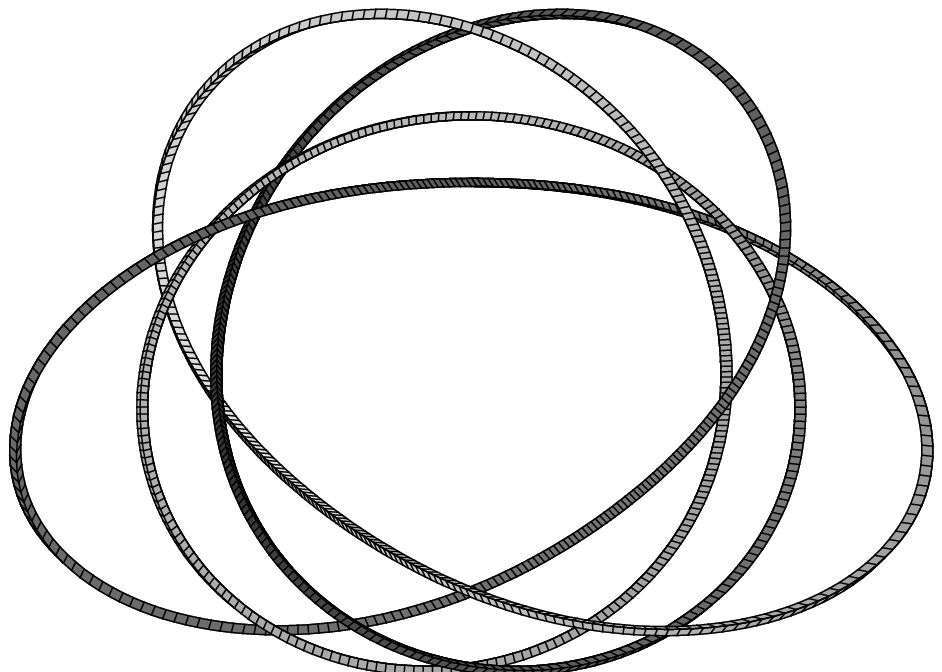}}
\medskip
\centerline{Figure $8$: The torus cable knot $(2,3)(2,3).$}
\endinsert

From the above we get the following partial description of the global
geometric monodromy. The typical regular fiber $F:=f^{-1}(c/2)$ is the
interior of the oriented surface obtained as the union of two pieces
$A$
and $B,$ where $A$ is a surface of genus $8$ with one boundary
component and
$B$ is a cylinder.
The pieces are glued together in the following way:
in each boundary component of $B$ there is an arc, which is glued to an
arc in the boundary of $A.$  The interior of $A$ or $B$ can be thought
of
as a Milnor fiber of the singularity at $0$ or $(-8,-4).$ So, the
geometric
monodromy around $0$ is a diffeomorphism with support in the interior
of
$A$ and $B,$ given for instance by a construction as in Paragraph $2.$
The piece $A$ can be constructed from the divide in Figure $9.$

\midinsert
\leftline{\epsffile{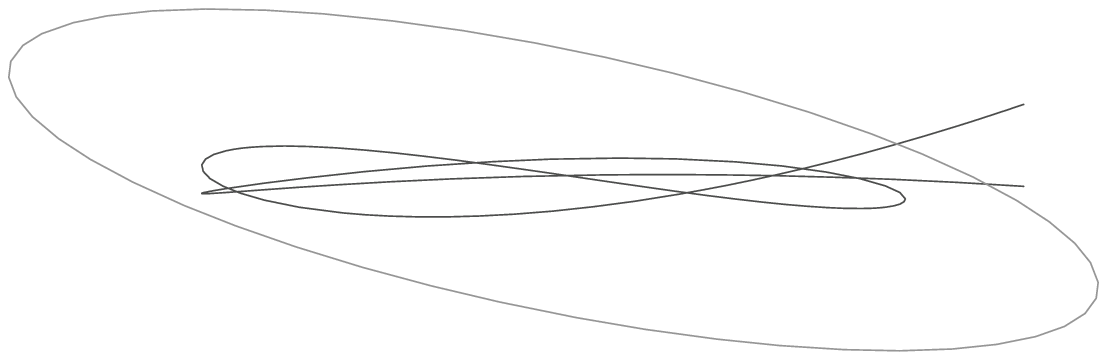}}
\medskip
\centerline{Figure $9$:  The curve $(x_s(t),y_s(t)),\, s:=1,$ as divide
for the singularity of $C.$}
\endinsert

Clearly, the monodromy in $B$ is a positive Dehn twist around
the simple essential closed curve $\delta_{17}$ in $B,$ whereas the
monodromy in $A$ is a product of positive Dehn twists around a system
$(\delta_1, \dots ,\delta_{16})$ of $16$ red, white or blue curves.
The monodromy around the critical point $c$ is a positive Dehn twist
around
a simple curve $\delta_{18},$ in $F,$ which is the union of two simple
arcs
$\alpha \subset A$ and $\beta \subset B.$ The arcs $\alpha$ and $\beta$
have
their endpoints $p,q \in A \cap B$ in common, and moreover the points
$p$ and
$q$ lie in different components of $A \cap B.$ The arc $\beta$ cuts
the curve
$\delta_{17}$ transversally in one point. The arc $\alpha$  intersects
the
curves $(\delta_1, \dots ,\delta_{16})$ transversally in some way.
For the position of the system $(\delta_{17},\, \beta)$ in B there is
up to
a diffeomorphism of the pair $(B,A \cap B)$ only one possibility. To
obtain
a complete description of the global monodromy  it remains to describe
the
position of the system $(\delta_1, \dots ,\delta_{16},\, \alpha)$
in $(A,A \cap B).$

We consider the family with parameter $s$ of parametrized curves with
para\-meter $t$:
$$x_s(t)=T(4,t)/8={t}^{4}-{t}^{2}+1/8,\,\,
y_s(t)=sT(6,t)/32+T(7,t)/64=$$
$$s{t}^{6}+{t}^{7}-3/2\,s{t}^{4}-7/4\,{t}^{5}+{{9}/{16}}\,s{t}^{2}+
{{7}/{8}}\,{t}^{3}-1/32\,s-{{7}/{64}}\,t,$$
where $T(d,t)$ is the Chebychev polynomial of degree $d.$ Let
$$
f_s(x,y):={s}^{4}{x}^{6}-{{3}/{128}}\,{s}^{4}{x}^{4}+{{1}/{1024}}\,{s}
^{4}{x}^{3}-2\,{s}^{2}{y}^{2}{x}^{3}-4\,sy{x}^{5}-
$$
$${x}^{7}+{{9}/{65536}}\,{s}^{4}{x}^{2}-{{3}/{262144}}\,{s}^{4}x+{{3}/{128}}
\,{s}^{2}{y}^{2}x-{{1}/{4096}}\,{s}^{2}{x}^{3}+
$$
$$
{{5}/{64}}\,s
y{x}^{3}+{{7}/{256}}\,{x}^{5}+{{1}/{4194304}}\,{s}^{4}-{{1}/{1024}}\,{s}^{2}{y}^{2}-{{1}/{1024}}\,sy{x}^{2}+
$$
$$
{y}^{4}+
{{3}/{1048576}}\,{s}^{2}x-
{{5}/{16384}}\,syx-{{7}/{
32768}}\,{x}^{3}-{{1}/{8388608}}\,{s}^{2}+
$$
$$
{{1}/{131072}}\,sy-{{1}/{4096}}\,{y}^{2}+{{7}/{16777216}}\,x+{{1}/{134217728}}
$$
be the equation, monic in $y,$ for the curve $(x_s(t),y_s(t),$ whose
real image is for $s=1$ a divide (see Figure $9$) for the singularity
of $C$ at $(0,0).$
The $0-$level of $f_s$ for $s=1$ consists of this divide and an
isolated
minimum not in one of its regions, which corresponds to the minimum of
the
restriction of $f$ to $\bR^2$ at $(-8,-4).$
For $a$ small, we call   $\delta_{17,a,s} \subset \{f_s=a\}$
the vanishing cycle of the local minimum of $f_s,s=1,$ which does not
belong
to a region. The curve $(x_s(t),y_s(t)), t \in \bR^2,s=7/24\,\sqrt{2},$
has
$8$ nodes, a cusp at $t=-1/2\,\sqrt{2},$(see Figure $10$) and at
infinity a singularity with
Milnor number 12. We now vary the parameter $s \in
[7/24\,\sqrt{2}-\sigma,1]$
from $1$ to $7/24\,\sqrt{2}-\sigma$ for a very small $\sigma >0.$

\midinsert
\cline{\epsffile{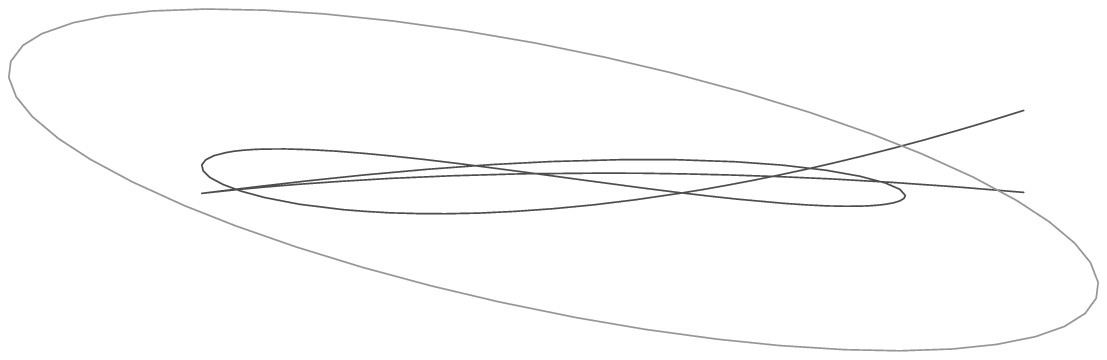}}
\medskip
\centerline{Figure $10$:  The curves $(x_s(t),y_s(t))$ for $
s:=7/24\,\sqrt{2}.$}
\endinsert

The value of the local minimum, which does not belong to a region,
becomes
smaller and by adjusting the parameter $a$ we can keep the cycle
$\delta_{17,a,s}$ in the new region which emerges from the cusp at
$s=7/24\,\sqrt{2}.$
Since the total Milnor number  of $f_s$ is
$(\hbox{\rm degree}(f_s)-1)(\hbox{\rm degree}(f_s)-2)-12=18,$
it follows that all its singularities have Milnor number $1$ for
$s=7/24\,\sqrt{2}-\sigma.$  The vanishing cycle of the node, which
appears when deforming the cusp singularity, will be called
$\delta_{18,s}$ and the vanishing cycle in the region of the divide of
Figure $9$, in whose boundary the cusp appeared, will be
called $\delta_{16,s}.$

\midinsert
\cline{\epsffile{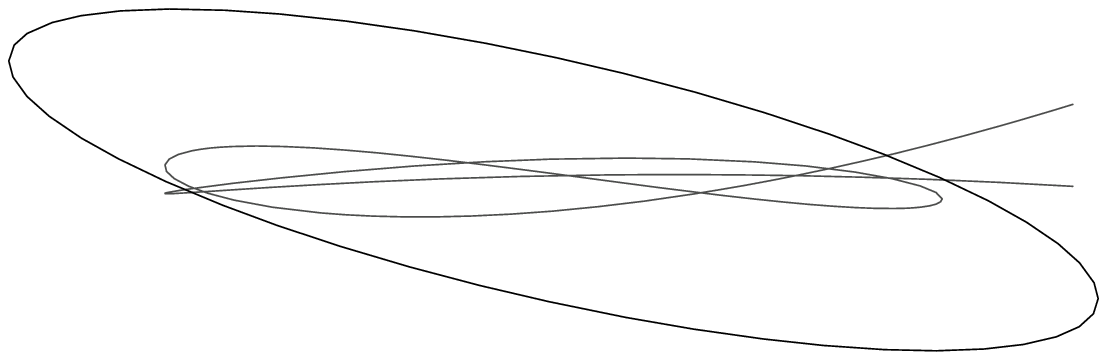}}
\medskip
\centerline{Figure $11$:  A divide, which does not come from a
singularity.}
\endinsert

We label the vanishing cycles  on the regular ribbon surface $F_{-}$ of
the divide of Figure $11$ by $\delta_{1}, ... ,\delta_{15}.$  The
cycles
$\delta_{17,a,s},\delta_{18,s},\delta_{16,s}$ deform without changing
their intersection pattern and $\delta_{16,s}$ becomes the cycle
$\delta_{16}$ of the regular ribbon surface $F_{+}$ of the divide of
Figure $9.$ Observe that the regular ribbon surface $F_{-}$ of the
divide
of Figure $11$ is naturally a subset of the regular ribbon surface
$F_{+}$ of the
divide of Figure $9.$ The description of the position of the system
$(\delta_1, \dots ,\delta_{16},\, \alpha)$ in $(A,A \cap B),$ for which
we are looking, is the system $(\delta_1, \dots ,\delta_{16})$ on
$F_{+},$ where the relative cycle $\alpha$ is a simple arc on
$F_{+}-F_{-}$ with endpoints on the boundary of $F_{+}$ and cutting
the cycle $\delta_{16}$ transversally in one point. Observe that
$F_{+}-F_{-}$ is a strip with core $\delta_{16}$ (see Figure $7$).

\goodbreak
\section{6}{Connected divides and fibered knots. Proof of Theorem 2}
In this section we assume, without loss of generality, that a divide
is linear and orthogonal near its crossing points. For a connected
divide $P \subset D(0,\rho),$ let $f_P:D(0,\rho) \to \bR$ be a generic
$C^{\infty}$ function,
such that $P$ is its $0$-level and that each region has exactly one
local maximum or minimum. Such a function exist for a connected divide
and is well defined up to sign and isotopy. In particular, there are no
critical points of saddle type other then the crossing points of the
divide. We assume moreover that the function $f_P$ is quadratic and
euclidean in a neighborhood of its critical points, i.e. for euclidean
coordinates (X,Y) with center at a critical point $c$ of $f_P$ we have
in a neighborhood of $c$ the expression $f_P(X,Y)=f_P(c)+XY$ or
$f_P(X,Y)=f_P(c)+X^2+Y^2.$ Let $\chi:D(0,\rho) \to [0,1]$ be a
$C^{\infty},$ positive function, which evaluates to zero outside of
the neighborhoods where $f_P$ is quadratic and to $1$ in some smaller
neighborhood of the critical points of $f_P.$ Let
$\theta_P:\partial B(0,\rho) \to \bC$ be given
by: $\theta_{P,\eta}(J(v)):=f_P(x)+
i \eta df_P(x)(u)-{1\over 2}\eta^2\chi(x)H_{f_P}(v)$ for
$J(v)=(x,u) \in TD(0,\rho)=D \times \bR^2$ and $\eta \in \bR, \eta >
0.$
Observe that the Hessian $H_{f_P}$ is locally
constant in a the neighborhood of the critical points of $f_P.$
The function $\theta_{P,\eta}$ is
$C^{\infty}$. Let
$\pi_{P,\eta}:\partial B(0,\rho) \setminus L(P) \to S^1$ be defined
by:$\pi_{P,\eta}(J(v)):=\theta_{P,\eta}(J(v))/
|\theta_{P,\eta}(J(v))|.$

\proclaim{Theorem 4}
Let $P \subset D(0,\rho)$ be a divide, such that the system of
immersed curves is connected.
The link $L(P)$ in
$\partial B(0,\rho)=\{J(v) \mid v \in T(D(0,\rho)) \ \text{ and
}\ \|J(v)\| =\rho\}$ is a fibered link. The map
$\pi_P:=\pi_{P,\eta}$ is for $\eta$ sufficiently small, a fibration
of the complement of $L(P)$ over $S^1.$  Moreover the fiber of the
fibration $\pi_P$ is $F(P)$ and the geometric monodromy is the product
of Dehn twist as in Theorem 3.
\endproclaim

The map $\pi_P$ is compatible with a regular product
tubular neighborhood of $L(P)$ in $\partial B(0,\rho).$ The map $\pi_P$
is a
submersion, so, since already a fibration near
$L(P),$ it is a fibration by a theorem of Ehresmann. The graphical
algorithm, see Figure $7$,  produces  in fact, up to a small isotopy of
the image, the projection of the
fiber $\pi_P^{-1}(i)$ on $D(0,\rho).$
This
projection is except above the twist of the strips a submersion.
The   proof of Theorem $4$ is given in the forthcoming paper [AC4] on
generic immersions of curves and knots.

\demo{ \bf Proof of Theorem $2$} The oriented fibered links $L(S)$ and $L(P)$
have the same geometric monodromies according to the Theorem 3 and 4. 
So, the links
$L(S)$ and $L(P)$ are
diffeomorphic.
\enddemo

\remark{\bf Remark}
Let $f(x,y)=0$ be a singularity $S,$ such that written in the
canonical coordinates of the charts of the embedded resolution
the branches of the
strict transform have equations of the form $u=a, a \in \bR.$ Let
$f_t(x,y), t \in [0,1]$ be a morsification, with its divide $P$ in
$D(0,\rho),$ obtained by blowing down generic real linear
translates of the strict transforms, as in [AC2]. We strongly
believe that with the use of [B-C1,B-C2] the following
transversallity property can be obtained,
and which we  state as a problem:\br
Their exists $\rho'_0>0,$ such that for all $t \in [0,1]$ and for all
$\rho' \in (0,\rho'_0]$ the $0$-levels of
$f_t(x,y)$ in $\bC^2$ meet transversally the boundary of
$$
B(0,\rho,\rho'):=\{(x+iu,y+iv) \in \bC^2 \mid x^2+y^2+u^2+v^2 \leq
\rho^2, u^2+v^2 \leq {\rho'}^2\}.
$$
It is easy to deduce from this transversallity statement an isotopy
between the links $L(S)$ and $L(P).$
\endremark

\remark{\bf Remark}
Bernard Perron has given a proof for the triviallity of the cobordism
from $L(S)$ to $L(P),$ which uses the holomorphic convexity of the
balls $B(0,\rho,\rho')$ of the previous remark [P].
\endremark

\par

\bigskip

\address{Universit\"at Basel \br Rheinsprung 21 \br CH-4051  Basel}

\topmatter
\title{Erratum: \br 
Real deformations and complex topology \br
of plane curve singularities}
\shorttitle{Real deformations and complex topology.}
\endtopmatter

In Section $5$ the parametrized curve $C$ should be $b(t):=(t^4,t^6+t^7)$
instead of $b(t):=(t^6+t^7,t^4)$ and accordingly $(-8,-4)$ has to
be $(-4,-8)$. We
intersect $C$ with the family of spheres
$S_r:=\{(x,y) \in \bC^2 \mid 4|x|^2+ |y|^2=r^2\}$.
For $0<r< 8\sqrt{2},$ the intersection
$K_r:=C \cap S_r$ is  the local knot in $S_r$ of the singularity at
$0 \in \bC^2$, at $r=8\sqrt{2}$ the knot $K_r$ is
singular with one
transversal crossing at $(-8,-4)$, and for $8\sqrt{2} < r$ the
knot $K_r$ is the so
called knot at infinity of the curve $C.$ Fig. $8$ of the text is a knot
projection of $K_r$ for small $r$. It is not possible to obtain
from this projection with only one crossing flip the type of
the knot $K_r$ for $r > 8\sqrt{2}$.
The figure  here below is the stereographic
knot projection of $K_r$ for $r = 8\sqrt{2}-1$,
which is not a minimal knot projection.
For $r=8\sqrt{2}$ the crossing at the bottom
flips and the knot $K_r, 8\sqrt{2} < r,$
becomes the $(4,7)$ torus knot.
The knot projection is a  braid projection, where
the axis is in the central pentagonal region.
The braid word is
$acabcaAabacabacabacab$ and flips at $r=8\sqrt{2}$ to
$acabcaaabacabacabacab$.

\midinsert
\cline{\epsffile{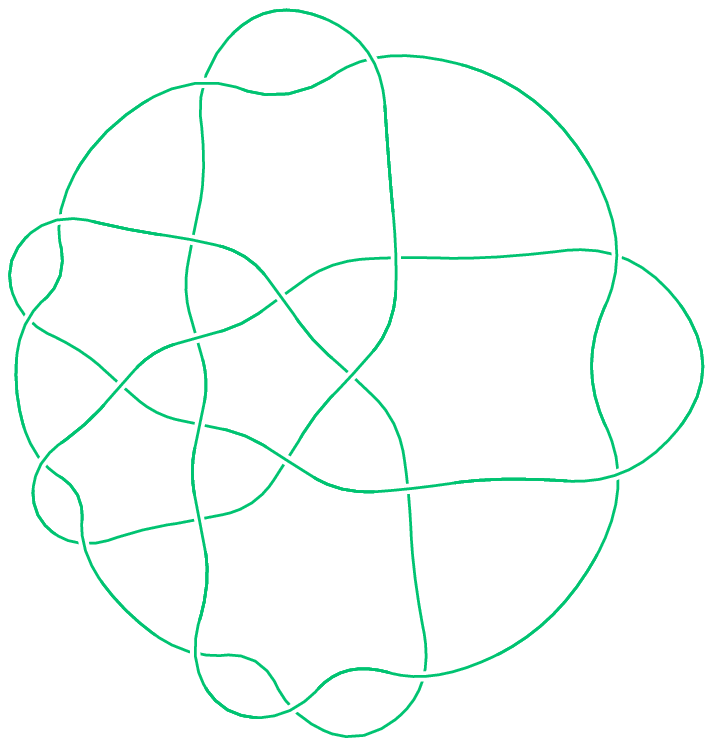}}
\medskip
\endinsert

\vskip 6.6cm  This picture was made with  KNOTSCAPE.

\bye